# THE STUDY OF CUCKOO OPTIMIZATION ALGORITHM FOR PRODUCTION PLANNING PROBLEM

A. Akbarzadeh[1], E. Shadkam[2]

[1,2]Department of Industrial Engineering, Faculty of Eng.; Khayyam University, Mashhad, Iran
[2]Department of Industrial Engineering, Isfahan University of Technology, Isfahan, Iran

## ABSTRACT

*Constrained Nonlinear programming problems are hard problems, and one of the most widely used and common problems for production planning problem to optimize. In this study, one of the mathematical models of production planning is survey and the problem solved by cuckoo algorithm. Cuckoo Algorithm is efficient method to solve continues non linear problem. Moreover, mentioned models of production planning solved with Genetic algorithm and Lingo software and the results will compared. The Cuckoo Algorithm is suitable choice for optimization in convergence of solution.*

## 1. INTRODUCTION

The optimization problems are so important in many different fields of science like physics, chemistry and engineering and their purpose is to find the best possible answer for an obvious problem.

The need of searching is the reason that many different searching algorithms have been provided. The evolutionary algorithms are a group of algorithms based on random searching that is inspired by the natural biological evolution modeling. They work on possible answers that have a superior feature and the longer generation survival that bring a closer estimate of the optimized answer. One of the newest evolutionary algorithms is the Cuckoo optimization algorithm. The current version of this algorithm is used to solve the continuous optimization problems and proved its performance in this field. Because of the ability of this algorithm in solving the continuous problems, it also has been discredited so that can solve the discrete issues too. The production planning is one of the problems that has a wide practical usage and is solvable in the discrete space. Therefore many expanded solving methods provide for this issue at all times. Considering the importance of this kind of problems, it comes to our mind to solve this problem using the Cuckoo algorithm.

One of the main models of production planning is to determine the lot sizing that has been investigated in different cases like multi period, multi item and that has capacity constraints with probable demand, removing the shortage cost and adding the level of services constraint by Bijari and at al. [1]

Also another model about determining the lot sizing and production scheduling with purpose of maximizing the profit has been investigated that in this model, the flexibility of choosing demand is imagined (the amount of chosen demands to comply is calculated by the model). In this article, a mathematical model has been suggested for PGLSP problem and it has been investigated [2].
Merzifonlou lu and Geunes studied the production planning problem with the hypothesis of flexibility in demand and the absence of capacity constraint. Not only they determined the amount of production in each period, but also they introduced a decision variable that can specify the percentage of accepted demand of each product in each period. The objective function of this problem is the same as the maximum profit objective function [6]. On the other hand, many algorithms and meta-heuristic algorithms have been provided in order to solve the production timing problem like "Ant colony optimization algorithm" [7], "Simulated annealing"[8], "Genetic algorithm" [9,10,11,12,13,14,15,16,17,18] and "Tabu search" [9]. Many articles are written about helpful nonlinear programming optimization [3], picture segmentations [4], Wind farms capacity determination [5], etc.

In the second section of this article, the Cuckoo algorithm will be explained, the third section belongs to analyzing the production planning model, the Cuckoo algorithm will be used in the mentioned problem in the fourth section and at last the conclusion will be provided.

## 2.1. Introducing the Cuckoo optimization algorithm

The Cuckoo optimization algorithm expanded by X. S. Yang and S. Deb in 2009 [20].

The Cuckoos laying method combined with the Levy Flight were the first clues of this algorithm creation instead of simple random isotropic hike. Later, R.Rajabioun investigated the Cuckoo optimization algorithm in detail in 2011 [21].

Like other evolutionary algorithms, this one begins with an initial population. This method works as the following steps. We assume that this population owns some eggs. First they put these eggs in other kind of birds' nests then wait until the host bird maintain these eggs beside her eggs. In fact, this lazy bird perfectly makes other birds to play an involuntary role in her generation survival. Some of the eggs that have less similarity to the host bird's eggs will be recognized and destroyed.

In fact the cuckoos improve and learn how to lay eggs more like the target host bird's eggs continuously and the host birds learn how to recognize the fake eggs.

More number of survived eggs in each zone shows more suitability of that zone, and more number of survived eggs, more attention pay to that zone and in fact this is the parameter that the Cuckoo optimization algorithm wants to optimize.

The Cuckoo optimization algorithm is as the rest:

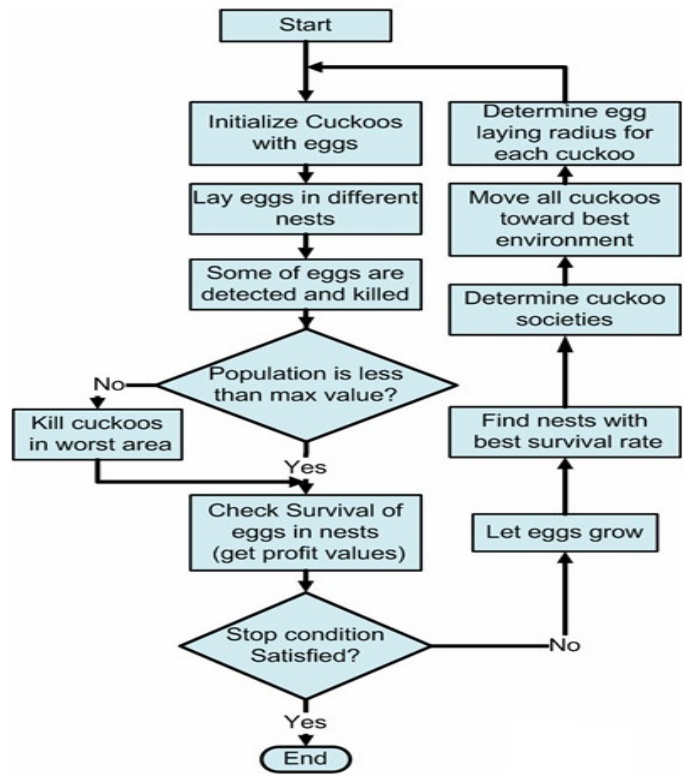

Figure 1: The Cuckoo optimization algorithm diagram

The variable values of the problem should take shape out of an array In order to solve an optimization problem. That array is called "habitat".

In an optimization problem, the next $N_{var}$ of a habitat will be a $1 \times N_{var}$ array that shows the current living location of cuckoos. This array describes as:

$$habitat = \left[ X_1, X_2, ..., X_{N\,var} \right]$$

The suitability (profit) in the current habitat obtains by evaluating the profit function ($f_p$) in the habitat. So:

$$profit = f\ b(habitat) = f\ b(X_1, X_2, ..., X_{N\,var})$$

A habitat matrix in size of $N_{pop} * N_{var}$ will be prepared for starting an optimization algorithm, then for each habitat a random number of eggs will be allocated.

The laying radius will be calculated by considering the number of eggs that each cuckoo lays and also the distance between the cuckoos and the current optimized zone. After that the cuckoos begin to lay in that zone. The laying radius calculates as:

$$ELR = a \times \frac{Number\ of\ current\ cuckoo's\ eggs}{Total\ number\ of\ eggs} \times (Var_{hi} - Var_{low})$$

Then each cuckoo begins to lay her eggs in the nests within her ELR

So after each laying round, p% of eggs (usually 10%) that is less profitable (their profit function is in lower level) destroys. Other chicks power up and grow in the host nests.

## 2.2. The cuckoo's migration

The cuckoos live in their environments while they are growing up and get older but when the laying time comes, they migrate to better habitats where the eggs have more chance to survive. After composing the groups in general different living locations (justified region or search space of the problem), the group with the best location will be targeted and other cuckoos will migrate there.

When the grown cuckoos live all around the environments, it is hard to find out each cuckoo belongs to which group. For solving this problem, the cuckoos will be grouped by the "K-means" method which is a classic way of grouping (finding a K between 3 and 5 is usually acceptable). When the cuckoos migrate to the target, they don't travel the direct way. They just travel (λ %) of the way with the deflection of (φ) as it is shown in figure 2.

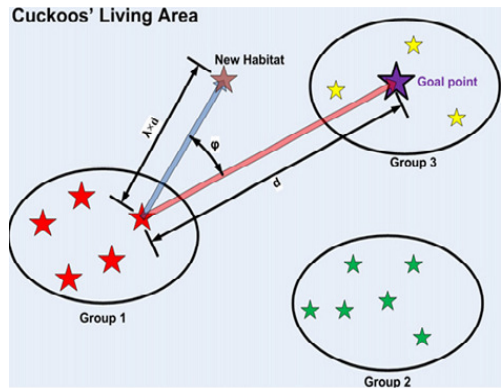

Figure 2: the cuckoos' migration to the target

These two parameters (λ,φ) helps cuckoos to explore a larger area. λ is a random number between 0 and 1 and φ is a number between $-\frac{\pi}{6}$ and $\frac{\pi}{6}$.

The migration formula is:

$$X_{Next\ Habitat} = X_{current\ Habitat} + F(X_{Goal\ Point} - X_{current\ Habitat})$$

After some iteration, all of the cuckoos gathered in an optimized point where the eggs have the most similarity to the host ones and access to the richest food sources is available. This location contains the most probable profit and the least possible number of killed eggs. Convergence more than 95% of all cuckoos in a single point ends the COA.

## 3. INTRODUCING THE MODEL AND ITS VARIABLES AND PARAMETERS

In this section, one of the production planning models will be introduced. The assumption is that a number of products have to be manufactured in a specific number of periods with the condition of minimizing its cost. Determining the optimized number of each product is the purpose. The parameters and the decision variables that is needed for the production planning problem, is given in the table 1 and 2.

Table1: the parameters of the model

| | |
|---|---|
| N | Number of all products |
| M | Number of all sources for producing |
| $P_{it}$ | Price of the product i in the period t |
| $C_{it}$ | Cost of the product i in the period t |
| $h_{it}$ | Maintenance cost of the product i in the period t |
| $Cap_{mt}$ | The source capacity m in the period t |
| $M_{it}$ | Maximum capacity of source i in the period t |
| $y_i$ | The occupied capacity for the product i |
| $D_{it}$ | The requirement of product i in the period t |
| $a_{im}$ | The intake of source m for the product i |
| $A_{it}$ | The preparing cost of the product i in the period t |

Table2: the decision variables of the model

| | |
|---|---|
| $s_{it}$ | The amount of sold product i in the period t |
| $x_{it}$ | The amount of product i in the period t |
| $B_{it}$ | Dearth of the product i in the period t |
| $I_{it}$ | Holding amount of the product i in the period t |
| $r_{it}$ | The delayed demand of the product i transported from the period t-1 to the period t |

The objective function purpose is maximizing the income of the sold products with reducing the production costs and maintenance and preparation costs. The objective function will be changed into a minimizing cost function by multiplying it to (-1), because many of the existing problems are minimizing cost problems and the model is for these kind of problems.

The mathematical model is at the rest:

$$\text{minimize } Z = -\sum_{i=1}^{N}\sum_{t=1}^{T}\text{pit.Sit} + \sum_{i=1}^{N}\sum_{t=1}^{T}\text{cit.X} + \sum_{i=1}^{N}\sum_{t=2}^{T}\text{hit.I} + \sum_{i=1}^{N}\sum_{t=1}^{T}\text{Zit.Ait} + \sum_{i=1}^{N}\sum_{t=1}^{T}\text{bit.Bit}$$

S.T

$$\sum_{i=1}^{N} y_i + I_{it} \leq M_{it}$$

$$\sum_{i=1}^{N} a_{im} . X_{it}$$

$$S_{it} = D_{it} + B_{it} + r_{it} . B_{i(t-1)}$$

$$X_{it} + B_{it} + I_{i(t-1)} - I_{it} - D_{it} - r_{it} . B_{i(t-1)} = 0$$

$$X_{it} \leq M . Z_{it}$$

In general, the first restriction is for the depository, the second one is for the used source, the third one is the selling balance, the fourth one is holding balance and the last one is the production with the condition of preparation.

We assume that the amount of demand for each product is apparent in each period and the delayed demand of each period will be transported to the next period. In this article, the production of 3 products in 5 periods is studied.

## 4. IMPLEMENTING THE CUCKOO ALGORITHM

There are number of input parameters in the Cuckoo optimization like other heuristic algorithm. The amounts of these parameters affect the final answer directly. The effects of these parameters are given in the table 3.

Table3: the changes of the input parameters

| Total Time | Cost | max Cuckoos | Cuckoos | Cluster |
|------------|------|-------------|---------|---------|
| 25.525 s | 6.5282E+06 | 30 | 7 | 3 |
| 12.591 s | 6.5099E+06 | 8 | 7 | 3 |
| 26.215 s | 6.5235E+06 | 8 | 7 | 6 |
| 28.172 s | 6.5257E+06 | 35 | 30 | 3 |
| 12.762 s | 6.5180E+06 | 8 | 7 | 2 |

So the parameters regulation is one of the most important factors that affect the final answer.

According to the amounts that are given in table 3 and our minimizing cost model, we calculated the best answer in the shortest possible time for the model. So the best input parameters for the model are given in table 4.

Table 4: the algorithm parameters

| Parameters | value |
|---|---|
| Number Of Initial Cuckoos | 7 |
| Minimum Number Of Eggs Laid by Each Cuckoo | 2 |
| Maximum Number Of Eggs Laid by Each Cuckoo | 4 |
| Maximum Number Of Iterations | 200 |
| Number Of Clusters for KNN | 3 |
| Maximum Number of Living Cuckoo | 8 |

Figure 3 shows the high isotropy of the algorithm in finding the answers of this problem.

This figure shows the high accuracy of the calculated answers of this algorithm. Comparing this algorithm with Genetic algorithm, the Cuckoo algorithm has fewer number of initial population and gives better answers in the shorter time.

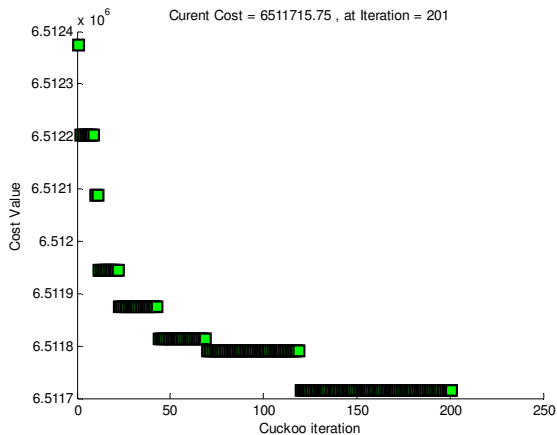

Figure3: the isotropy speed of the algorithm for 200 iterations

If checking the operation of meta-heuristic algorithms for finding the optimized answer and determining the production set would be in order, the model should be solved with the default parameters with different methods. The general conclusion will be obtained by comparing the results.

Also it is so obvious that if we implement a single meta-heuristic algorithm for a single problem a few times, we probably will come up with different answers for the optimized (or near optimized) quantity of the objective function. It is because of many meta-heuristic algorithms such as Genetic and the Cuckoo algorithm use the random operators in order to find the optimized answer.

Because the model is nonlinear and doesn't have the optimized answer, we solved this problem by Genetic algorithm and the Lingo 11 software. After 20 iterations, the average of the answers is given in the table 5.

Table5: comparing the average of the answers after 20 iterations

| GA | COA | | | Lingo |
|---|---|---|---|---|
| average | deviation | best | average | average |
| 2.38e+11 | 1.38e+04 | 6.50e+06 | 6.52e+06 | 2.69e+09 |
| 18 s | ------- | 16.674 s | 16.935 s | 1 s |

## 5. CONCLUSION

In this article, two of the meta-heuristic algorithms named Genetic algorithm and the Cuckoo algorithm have been used for solving one of the production planning problems. Because the Cuckoo algorithm is the expanded version of the Genetic algorithm, it gives better answers.
Also the Lingo software provides an answer near the optimized answer using the Global Solver ability but with increasing the scale of the problem, the calculation timing will be incomputable.
The Genetic algorithm cannot find the optimized answer but it finds the answer in an acceptable duration. The Cuckoo algorithm finds better answers of this problem in the suitable duration but we cannot name it as the best way finding answer in contrast of other algorithms because each year, many algorithms are published that covers the defects of earlier algorithms. But in this problem, the Cuckoo algorithm could easily find the answer with lower number of initial population and less iteration.

The subject of this article can be expanded in many other ways such as solving the model with other algorithms or adding the probable demand and the transportation cost restriction to the model in case of comparing two algorithms.